\documentclass[a4paper,11pt]{article}

\usepackage[latin1]{inputenc}
\usepackage{amsfonts,amsmath,amssymb,amsthm,graphicx,epsfig,float,xcolor}
\usepackage[T1]{fontenc}
\usepackage[english,francais]{babel}
\usepackage[all]{xy}

\usepackage{amsthm}
\newtheorem{Thm}{Theorem}

\newtheorem{Lem}{Lemma}
\newtheorem{Cor}{Corollary}

\theoremstyle{definition}
\newtheorem{Def}{Definition}
\theoremstyle{remark}
\newtheorem{Rem}{Remark}
\newtheorem*{Proof}{Proof}

\newcommand{\Cc}{\mathbb{C}}
\newcommand{\Nn}{\mathbb{N}}

\newcommand{\Pp}{\mathbb{P}}

\title{{\bf The bifurcation measure is exponentially mixing}}
\author{Henry De Th\'elin}
\date{}

\begin{document}
\maketitle

\def\figurename{{Fig.}}%
\def\proofname{Preuve}
\def\contentsname{Sommaire}%

\selectlanguage{english}
\begin{center}
{\bf{ }}
\end{center}

\begin{abstract}

We prove general mixing theorems for sequences of meromorphic maps on compact K\"ahler manifolds. We deduce that the bifurcation measure is exponentially mixing for a family of rational maps of $\Pp^q(\Cc)$ endowed with suitably many marked points.

\end{abstract}



Classification: 37A25, 37F46, 37F80.

\section*{{\bf Introduction}}
\par

Given a holomorphic endomorphism of $\Pp^k(\Cc)$, $f$, of degree $d \geq 2$, Forn{\ae}ss and Sibony defined the Green current $T$ associated with $f$ (see \cite{FS} and \cite{FS1}), whose support is the Julia set of $f$, that is, the set of points $x \in \Pp^k(\Cc)$  for which the sequence $(f^n(x))_n$ is not normal in some neighborhood of $x$. This current has a continuous potential: we can therefore define its self-intersection $\mu=T^k$ (see \cite{FS}). The measure $\mu$ obtained in this way is mixing (see \cite{FS}), it is the unique measure of maximal entropy $k \log (d)$ (see \cite{BD2}), and its Lyapunov exponents are bounded from below by $\frac{\log(d)}{2}$ (see \cite{BD1}).

In a similar way, let now $\Lambda$ be a complex K\"ahler manifold and $\widehat{f}: \Lambda \times \Pp^1 \longrightarrow \Lambda \times \Pp^1$ an algebraic family of rational maps of degree $d \geq 2$: $\widehat{f}$ is holomorphic and $\widehat{f}(\lambda,z)=(\lambda, f_{\lambda}(z))$ where $f_{\lambda}$ is a rational map of degree $d$. Let $a$ be a marked point, i.e., a rational function $a: \Lambda \longrightarrow \Pp^1$. As for holomorphic endomorphisms, a fundamental question is to study the bifurcation locus, that is, the set of parameters $\lambda_0 \in \Lambda$ for which the sequence $(f^n_{\lambda}(a(\lambda)))_n$ is not normal in some neighborhood of $\lambda_0$.
For example, in the historical example, $f_{\lambda}(z)=z^d + \lambda$ with $\lambda \in \Cc$ and $a(\lambda)=\lambda$, the bifurcation locus is the Mandelbrot set.

DeMarco introduced in \cite{DeM} a current of bifurcation $T_{\rm bif}$ on $\Lambda$: it is a positive
closed current of bidegree $(1, 1)$ whose support is exactly the bifurcation locus and Bassanelli and Berteloot (\cite{BB1}) then defined its self-intersections $T^l_{\rm bif}$. The maximal intersection $\mu_{\rm bif} := T_{\rm bif}^{\rm dim(\Lambda)}$ is known as the bifurcation measure and in the mentioned historical example, it corresponds to the harmonic measure of the Mandelbrot set.
For this harmonic measure, Graczyk and \'Swiatek (see \cite{GS}) proved that the Lyapounov exponent of $\mu_{\rm bif}$ is equal to $\log(d)$. In \cite{DGV2}, we partially extended this result to the case of any pair $(f,a)$ and with a quasi-projective variety $\Lambda$ of dimension $1$ by showing that the Lyapounov exponent is bounded from below by $\frac{\log(d)}{2}$. Moreover, we extended the notion of entropy to the context of more general parameter families in \cite{DGV}, and proved that the measure  $\mu_{\rm bif}$ has maximal entropy.

In this article, we continue the analogy with endomorphisms, showing that the measure $\mu_{\rm bif}$ is mixing in a very general setting. In \cite{GKN}, Ghioca, Krieger and Nguyen  proved that the Mandelbrot set is not the Julia set of a polynomial map (also refer to \cite{L} for another similar result). So let us clarify the context and what we mean by "mixing".

Let $\Lambda$ be a smooth complex quasi-projective variety and let $\widehat{f}: \Lambda \times \Pp^q \longrightarrow \Lambda \times \Pp^q$ be an algebraic family of endomorphisms of $\Pp^q$ of degree $d \geq 2$: $\widehat{f}$ is a morphism and $\widehat{f}(\lambda,z)=(\lambda, f_{\lambda}(z))$ where $f_{\lambda}$ is an endomorphism of $\Pp^q$ of algebraic degree $d$.

Assume that the family $\widehat{f}$ is endowed with $k$ marked points $a_1, \dots , a_k : \Lambda \longrightarrow \Pp^q$ where the $a_i$ are morphisms and  suppose that ${\rm dim}( \Lambda)=qk$. As in the case $k=q=1$, we can define $T_{\rm bif}$ and $\mu_{\rm bif}=T_{\rm bif}^{kq}$ (see the paragraph \ref{background} for more details, or refer to \cite{DGV} and \cite{Du}). 

For $n \in \Nn$, define 
$$\mathfrak{a}_n(\lambda)=(f_{\lambda}^n(a_1(\lambda)), \cdots , f_{\lambda}^n(a_k(\lambda))) \mbox{    }   \mbox{   ,   } \mbox{    } \lambda \in \Lambda.$$

Let $\iota: \Lambda \hookrightarrow \Pp^m$ be an embedding of $\Lambda$ into a complex projective space. We identify $\Lambda$ with $\iota(\Lambda)$ and we denote by $\overline{\Lambda}$ the closure of $\Lambda$ in $\Pp^m$. Let $p_n: \overline{\Lambda} \longrightarrow (\Pp^q)^k $ be the meromorphic map obtained by taking the closure of the graph of $\mathfrak{a}_n$ in $\overline{\Lambda} \times (\Pp^q)^k$.

In this context, we first have (see paragraph \ref{preliminaries} for a review of notions used here, such as dsh functions, locally moderate measures and PB probability measures) 

\begin{Thm}{\label{Th2}}

We assume that $\mu_{\rm bif} \neq 0$ and ${\rm dim}( \Lambda)=qk$. Take $U \subset (\mathbb{P}^q)^k$ an open set and $\mu$ a locally moderate positive measure in $U$.

Let $V \Subset U$ and $W \Subset \Lambda$ be relatively compact sets. Then for $s \in ]1, + \infty[$ and $0 \leq \nu \leq 2$, there exists a positive constant $C$ such that

$$\left| \int \psi(p_n) \varphi \frac{p_n^*( \mu_{| \overline{V}})}{d^{kqn}} - \int \psi d \mu_{| \overline{V}} \int \varphi  d \mu_{\rm bif}  \right| \leq Cd^{-n \nu/2}   | \psi |_{L^s(\mu_{| \overline{V}})} | \varphi |_{C^{\nu}}$$

for every $n \in \Nn$, $\psi \in DSH((\mathbb{P}^q)^k)$ and $\varphi \in C^{\nu}$ with compact support in $\overline{W}$. 

\end{Thm}

We have $\mu_{\rm bif}= T_{\rm bif}^{kq}$ where the $T_{\rm bif}$ has locally H\"older potentials, hence the measure $\mu_{\rm bif}$ is locally moderate (see \cite{DSN}), and we deduce

\begin{Cor}

We assume that $\mu_{\rm bif} \neq 0$ and that $\Lambda$ is a Zariski open set in $(\mathbb{P}^q)^k$. 

Let $V \Subset \Lambda=U$ be a relatively compact set. Then for $s \in ]1, + \infty[$ and $0 \leq \nu \leq 2$, there exists a positive constant $C$ such that

$$\left| \int \psi(p_n) \varphi \frac{p_n^*( {\mu_{\rm bif}}_{| \overline{V}})}{d^{kqn}} - \int \psi d {\mu_{\rm bif}}_{| \overline{V}} \int \varphi  d {\mu_{\rm bif}}_{| \overline{V}}  \right| \leq Cd^{-n \nu/2}   | \psi |_{L^s({\mu_{\rm bif}}_{| \overline{V}})} | \varphi |_{C^{\nu}}$$

for every $n \in \Nn$, $\psi \in DSH((\mathbb{P}^q)^k)$ and $\varphi \in C^{\nu}$ with compact support in $\overline{V}$. 

\end{Cor}

This means that the bifurcation measure is exponentially mixing.

To prove the Theorem \ref{Th2}, we first establish a very general mixing theorem for sequences of meromorphic maps. The proof follows the approach used by Dinh, Nguyen and Sibony to prove stochastic properties for holomorphic endomorphisms in $\Pp^k(\Cc)$ (see \cite{DSN} and also \cite{DS08}). Let us explain the context.

Let $(X,\omega)$ and $(X',\omega')$ be compact K\"ahler manifolds of dimension $l$ and consider a sequence of dominant meromorphic maps $p_n:X \longrightarrow X'$. 

Fix $U \subset X'$ an open set, and let $\mu$ be a locally moderate positive measure on $U$.

Then, we have

\begin{Thm}{\label{Th1}}

Let $V \Subset U$ be a relatively compact set and take $\lambda$ a PB probability measure on $X'$ . Then for $s \in ]1, + \infty[$ and $0 \leq \nu \leq 2$, there exists constants $C_1 , C_2$ such that

$$\left| \int \psi(p_n) \varphi \frac{p_n^*( \mu_{| \overline{V}})}{\delta_l(p_n)} - \int \psi d \mu_{| \overline{V}} \int \varphi \frac{p_n^*( \lambda) }{\delta_l(p_n)} \right| \leq C_1 \frac{\delta_{l-1}(p_n)}{\delta_l(p_n)} | \psi |_{L^s(\mu_{| \overline{V}})} | \varphi |_{DSH}$$

for every  $n \in \Nn$, $\psi \in DSH(X')$ and $\varphi \in DSH(X)$ and

$$\left| \int \psi(p_n) \varphi \frac{p_n^*( \mu_{| \overline{V}})}{\delta_l(p_n)} - \int \psi d \mu_{| \overline{V}} \int \varphi \frac{p_n^*( \lambda) }{\delta_l(p_n)} \right| \leq C_2 \left( \frac{\delta_{l-1}(p_n)}{\delta_l(p_n)} \right)^{\nu/2} | \psi |_{L^s(\mu_{| \overline{V}})} | \varphi |_{C^{\nu}}$$

for every $n \in \Nn$, $\psi \in DSH(X')$ and $\varphi \in C^{\nu}(X)$.

Here $\delta_l(p_n)=\int p_n^*(\omega'^l)$ and $\delta_{l-1}(p_n)=\int p_n^*(\omega'^{l-1}) \wedge \omega$.

\end{Thm}

Here is the outline of this paper: In the first paragraph, we review the various notions used in these statements (locally moderate measures, PB probability measures, etc.) and explain why the integrals in Theorem \ref{Th1} are well-defined. Then we will demonstrate Theorem \ref{Th1}. In the second paragraph, we begin with some background about parameter families and we will prove Theorem \ref{Th2}, using crucially Theorem \ref{Th1} and pluripotential theory.

{\bf Acknowledgments:} Thanks to Gabriel Vigny for many useful discussions on this article.

\section{\bf Proof of Theorem \ref{Th1}}{\label{proof-Th1}}

\subsection{\bf Preliminaries}{\label{preliminaries}}

Let $(X, \omega)$ be a compact  K\"ahler manifold. We start with some reminders about dsh functions, PB measures and locally moderate measures  (see \cite{DS08}).

A function $\varphi$ is quasi-plurisubharmonic (qpsh) if it is locally written as
the sum of a psh function and a $C^{\infty}$ one. Such a function verifies $dd^c \varphi \geq -c \omega$ in the sense of currents for a constant $c \geq 0$. A set of $X$ is said to be pluripolar if it is contained in $\{ \varphi = - \infty \}$ where $\varphi$ is
a qpsh function. We call dsh function, any function defined outside a pluripolar set, which is written as the difference of two qpsh functions. Let us denote DSH(X) the set of dsh functions on $X$. If $\varphi$ is a dsh function, there are two positive closed currents $T^{\pm}$ of bidegree $(1, 1)$ such that $dd^c \varphi = T^+ - T^-$. We can then define a
norm (see \cite{DSN} paragraph $3$):

$$\| \varphi \|_{\mbox{DSH}} := \| \varphi \|_{L^1(X)} + \inf \| T^{\pm} \|$$

with $T^{\pm}$ as above.

A positive measure $\mu$ is PB if qpsh functions are integrable with respect to this measure. In particular, PB measures have no mass on pluripolar sets. Let $\mu$ be a non-zero PB positive measure on $X$. For $\varphi \in DSH(X)$, define

$$\| \varphi \|_{\mu} := | \langle \mu, \varphi \rangle | + \inf \| T^{\pm} \|$$

with $T^{\pm}$ as above.

The semi-norm $\| . \|_{\mu}$ is in fact a norm on DSH(X) which is equivalent to $\| . \|_{\mbox{DSH}}$ (see Proposition A.4.4 in \cite{DS08}).

The measure $\mu$ is said to be locally moderate (see \cite{DSN}) if for any open set $U \subset X$, any compact set $K \subset U$ and any compact family $\mathcal{F}$ of psh functions on $U$, there are constants $\alpha >0$ and $c>0$ such that

$$\int_K e^{- \alpha \varphi} d \mu \leq c  \mbox{   }   \mbox{   }  \mbox{  for  }  \mbox{   } \varphi \in \mathcal{F}.$$

By using Proposition 2.1 in \cite{DS06} and Cauchy-Schwarz inequality, if $\mu$ is locally moderate, for any open set $U \subset X$, any compact set $K \subset U$ and any compact family $\mathcal{D}$ of dsh functions on $X$, there are constants $\alpha >0$ and $c>0$ such that

$$\int_K e^{ \alpha | \varphi |} d \mu \leq c  \mbox{   }   \mbox{   }  \mbox{  for  }  \mbox{   } \varphi \in \mathcal{D}.$$

\subsection{\bf About the definition of the integrals in Theorem \ref{Th1}}

We begin by showing that all the integrals in Theorem \ref{Th1} are well-defined.

Take $\xi \leq 0$ a qpsh function and $K$ a compact set. Then, since $\mu$ is locally moderate, there exists $\alpha>0$ such that $\int_K -\alpha \xi d \mu \leq \int_K e^{- \alpha \xi} d \mu < + \infty$. It implies that $\mu$ gives no mass to analytic subsets in $X'$ and integrates dsh functions ($\mu$ is PB).

Let $\Gamma_{p_n}$ be the graph of $p_n$ in $X \times X'$ and $\alpha_1, \alpha_2$ the projections of $\Gamma_{p_n}$ onto $X$ and $X'$, respectively. Take $\varphi$ a continuous map, then ${\alpha_2}_*( \alpha_1^*  \varphi) $ is continuous outside an analytic subset of $X'$. Hence, if $\lambda$ is a PB positive measure on $X'$, as it gives no mass to analytic subset in $X'$, we can define a positive measure $p_n^* \lambda$ with the formula:

$$\langle p_n^* \lambda , \varphi \rangle:=\langle \lambda, {\alpha_2}_*( \alpha_1^*  \varphi )\rangle.$$

Now, let $\varphi$ be a DSH function in $X$. With the above notations, ${p_n}_* \varphi= {\alpha_2}_*( \alpha_1^*  \varphi)$ is DSH (see paragraphs 2.3 and 2.4 in \cite{DS06}). By definition, since $\lambda$ is a PB probability measure, qpsh functions are $\lambda$-integrable, hence $\langle  p_n^* \lambda ,  \varphi \rangle:=  \langle   \lambda , {p_n}_* \varphi \rangle$ is well-defined.

Finally, let $\psi_1,\psi_2$ be two negative qpsh functions in $X'$. We have

\begin{equation*}
\begin{split}
\int \left| \psi_1  \psi_2 \right|  d \mu_{| \overline{V}} &\leq \left( \int |\psi_1 |^2  d \mu_{| \overline{V}} \right)^{1/2} \left( \int  |\psi_2 |^2 d \mu_{| \overline{V}} \right)^{1/2}\\
&\leq \frac{2}{\alpha^2} \left( \int e^{\alpha |\psi_1 |}d \mu_{| \overline{V}} \right)^{1/2}  \left( \int e^{\alpha |\psi_2 |}d \mu_{| \overline{V}} \right)^{1/2} < + \infty\\
\end{split}
\end{equation*}

where we used that $\mu$ is locally moderate and the inequality $\frac{\alpha^2 x^2}{2} \leq e^{\alpha x}$ for $x \geq 0$.

We deduce that $\int \psi(p_n) \varphi \frac{p_n^*( \mu_{| \overline{V}})}{\delta_l(p_n)} := \int \psi   \frac{{p_n}_*(\varphi)}{\delta_l(p_n)}  d \mu_{| \overline{V}}$ is well defined for $\psi \in DSH(X')$ and $\varphi \in DSH(X)$, since ${p_n}_* \varphi$ is DSH in $X'$.

\subsection{\bf Proof of Theorem \ref{Th1}:}

In this paragraph, we prove Theorem \ref{Th1}. We follow the ideas of Dinh-Sibony-Nguyen used in \cite{DSN} to prove that the measure of maximal entropy for a holomorphic endomorphism of $\Pp^k(\Cc)$ is exponentially mixing. In particular, in what follows, we use their notation $\Lambda_n(\varphi):= \frac{{p_n}_*(\varphi)}{\delta_l(p_n)}$. For $\psi \in DSH(X')$ and $\varphi \in DSH(X)$, we have

\begin{equation*}
\begin{split}
\left| \int \psi(p_n) \varphi \frac{p_n^*( \mu_{| \overline{V}})}{\delta_l(p_n)} - \int \psi d \mu_{| \overline{V}} \int \varphi \frac{p_n^*( \lambda) }{\delta_l(p_n)} \right| &=\left| \int \psi \Lambda_n(\varphi) d \mu_{| \overline{V}} - \int \psi d \mu_{| \overline{V}} \int \Lambda_n(\varphi) d \lambda \right| \\
&=\left| \int \psi \left(\Lambda_n(\varphi) - \int \Lambda_n(\varphi) d \lambda \right) d \mu_{| \overline{V}}  \right| \\
&=\left| \int \psi \widetilde{\Lambda_n}(\varphi)  d \mu_{| \overline{V}}  \right| \\
\end{split}
\end{equation*}

where we denote $ \widetilde{\Lambda_n}(\varphi)= \Lambda_n(\varphi) - \langle \Lambda_n(\varphi) , \lambda \rangle$. We can estimate the norm of this function by using the following Lemma.

\begin{Lem}

There exists a constant $r>0$ which depends only on $X$ such that

$$\left| \widetilde{\Lambda_n}(\varphi) \right|_{DSH(\lambda)} \leq r |\varphi|_{DSH} \frac{\delta_{l-1}(p_n)}{\delta_l(p_n)}$$

for every $\varphi \in DSH(X)$.

\end{Lem}

\begin{Proof}

Write $dd^c \varphi= S^+ - S^-$ with $S^{\pm}$ positive closed $(1,1)$-currents. We have

$$dd^c \widetilde{\Lambda_n}(\varphi)= dd^c \Lambda_n(\varphi)=\frac{{p_n}_* (dd^c \varphi)}{\delta_l(p_n)}= \frac{{p_n}_* (S^{+})}{\delta_l(p_n)} -\frac{{p_n}_* (S^{-})}{\delta_l(p_n)}.$$

In the above, the push-forward ${p_n}_*(S^{\pm})$ is well-defined using Meo's result (see \cite{M}) and the projections from the graph of $p_n$ (taking a desingularization if necessary), since $S^{\pm}$ are positive $(1,1)$-currents.

By definition

$$\left| \widetilde{\Lambda_n}(\varphi) \right|_{DSH(\lambda)} = \left| \langle \widetilde{\Lambda_n}(\varphi), \lambda \rangle  \right| + \min \| R^{\pm} \|$$

where the minimum is taken on  positive closed $(1,1)$-currents $R^{\pm}$ such that $dd^c \widetilde{\Lambda_n}(\varphi)= R^+ - R^-$. Hence,

$$\left| \widetilde{\Lambda_n}(\varphi) \right|_{DSH(\lambda)} \leq \left\| \frac{{p_n}_* (S^{\pm})}{\delta_l(p_n)} \right\|$$

and it remains to estimate the norm of the term on the right.

By using Proposition $2.2$ in \cite{DS06}, there exists a constant $r>0$ that depends only on $X$ such that $S^+= \beta + dd^c u$ with $\beta$ a  smooth $(1,1)$-form, $u$ a qpsh function and 

$$-r  \| S^+ \|  \omega \leq  \beta \leq  r  \| S^+ \|  \omega.$$

Finally,

$$\left\| \frac{{p_n}_* (S^{+})}{\delta_l(p_n)} \right\|= \langle \frac{{p_n}_* (S^{+})}{\delta_l(p_n)} , \omega'^{l-1} \rangle=    \langle \beta,  \frac{{p_n}^* (\omega'^{l-1} )}{\delta_l(p_n)}  \rangle \leq r \|S^+\| \frac{\delta_{l-1}(p_n)}{\delta_l(p_n)}$$

and the Lemma follows.

\end{Proof}

We continue the proof of Theorem \ref{Th1}.

Take $r>0$ such that $\frac{1}{s}+ \frac{1}{r}=1$. By using Proposition A.4.4 in \cite{DS08}, the sequence

$$\left( \widetilde{\Lambda_n}(\varphi) \frac{\delta_{l}(p_n)}{\delta_{l-1}(p_n)} \right)=\left( (\Lambda_n(\varphi) - \langle \Lambda_n(\varphi) , \lambda \rangle)\frac{\delta_{l}(p_n)}{\delta_{l-1}(p_n)} \right)$$

is bounded in $DSH(X)$ for $ \varphi \in DSH(X)$ with $|\varphi |_{DSH} \leq 1$ (recall that $\lambda$ is a PB probability measure). Hence, since $\mu$ is locally moderate and $\overline{V}$ is a compact set, there exist positive constants $\alpha , C$ such that $\langle e^{\alpha |\psi|} , \mu_{| \overline{V}} \rangle \leq C$ for every $\psi$ in the above sequence.

By using the inequalities $x^r \leq r! e^x \leq r^r e^x$ for $x \geq 0$ (consider integer parts of $r$, if necessary), we obtain 

$$\left\langle \left( \alpha | \Lambda_n(\varphi) - \langle \Lambda_n(\varphi) , \lambda \rangle|  \frac{\delta_{l}(p_n)}{\delta_{l-1}(p_n)} \right)^r , \mu_{| \overline{V}} \right\rangle  \leq C r^r$$

for $ \varphi \in DSH(X)$ with $|\varphi |_{DSH} \leq 1$. Thus, 

\begin{equation}{\label{eq1}}
 | \Lambda_n(\varphi) - \langle \Lambda_n(\varphi) , \lambda \rangle|_{L^r(\mu_{| \overline{V}})} \leq \frac{rC^{1/r}}{\alpha}  \frac{\delta_{l-1}(p_n)}{\delta_{l}(p_n)} |\varphi |_{DSH}
\end{equation}

for every $ \varphi \in DSH(X)$.

We deduce,

\begin{equation*}
\begin{split}
\left| \int \psi(p_n) \varphi \frac{p_n^*( \mu_{| \overline{V}})}{\delta_l(p_n)} - \int \psi d \mu_{| \overline{V}} \int \varphi \frac{p_n^*( \lambda) }{\delta_l(p_n)} \right| &= \left| \int \psi \left(\Lambda_n(\varphi) - \int \Lambda_n(\varphi) d \lambda \right) d \mu_{| \overline{V}}  \right| \\
& \leq | \psi |_{L^s(\mu_{| \overline{V}})} | \Lambda_n(\varphi) - \langle \Lambda_n(\varphi) , \lambda \rangle|_{L^r(\mu_{| \overline{V}})} \\
& \leq \frac{rC^{1/r}}{\alpha}  \frac{\delta_{l-1}(p_n)}{\delta_{l}(p_n)}  | \psi |_{L^s(\mu_{| \overline{V}})}   |\varphi |_{DSH}\\
\end{split}
\end{equation*}

which proves the first inequality.

\begin{Rem}

As mentioned above, if $\psi$ is DSH with $| \psi |_{DSH} \leq 1$, we have $| \psi |^s \leq \left( \frac{s}{\alpha} \right)^s e^{\alpha | \psi|} $ and $\langle e^{\alpha |\psi|} , \mu_{| \overline{V}} \rangle \leq C$. Hence, there exists a positive constant $C'$ which does not depend on $\psi$ such that $| \psi |_{L^s(\mu_{| \overline{V}})} \leq C' | \psi |_{DSH}$. It means that we can replace $| \psi |_{L^s(\mu_{| \overline{V}})} $ with $| \psi |_{DSH}$ in the first inequality of Theorem \ref{Th1}.

\end{Rem}{\label{rem1}}

The second inequality follows classically from the theory of interpolation between the Banach spaces $\mathcal{C}^0$ and  $\mathcal{C}^2$ (see \cite{DS08} p.34).

Let $L: \mathcal{C}^0 \longrightarrow  L^r(\mu_{| \overline{V}})$ be the linear operator $L( \varphi)= \Lambda_n(\varphi) - \langle  \Lambda_n(\varphi), \lambda \rangle$. We have

$$| \Lambda_n(\varphi) - \langle  \Lambda_n(\varphi), \lambda \rangle |_{L^r(\mu_{| \overline{V}})} \leq 2 | \varphi |_{\mathcal{C}^0} \mu(\overline{V})^{1/r}$$

for $\varphi \in \mathcal{C}^0$. Since $| \varphi |_{DSH} \lesssim | \varphi |_{\mathcal{C}^2} $ for $\varphi \in \mathcal{C}^2$, by using the inequality (\ref{eq1}), we deduce

\begin{equation*}
\begin{split}
 | \Lambda_n(\varphi) - \langle \Lambda_n(\varphi) , \lambda \rangle|_{L^r(\mu_{| \overline{V}})} &\leq \frac{rC^{1/r}}{\alpha}  \frac{\delta_{l-1}(p_n)}{\delta_{l}(p_n)} |\varphi |_{DSH}\\
 & \leq C''  \frac{\delta_{l-1}(p_n)}{\delta_{l}(p_n)} |\varphi |_{\mathcal{C}^2}.\\
 \end{split}
\end{equation*}

By applying the theory of interpolation to $L$, for $0 \leq \nu \leq 2$, there exists a positive constant $A_{\nu}$, independent of $L$ such that

$$ | \Lambda_n(\varphi) - \langle \Lambda_n(\varphi) , \lambda \rangle|_{L^r(\mu_{| \overline{V}})} \leq A_{\nu} \left( 2  \mu(\overline{V})^{1/r} \right)^{1- \nu/2} \left(  C''  \frac{\delta_{l-1}(p_n)}{\delta_{l}(p_n)}  \right)^{\nu /2} |\varphi |_{\mathcal{C}^{\nu}}.$$

We obtain

\begin{equation*}
\begin{split}
&\left| \int \psi(p_n) \varphi \frac{p_n^*( \mu_{| \overline{V}})}{\delta_l(p_n)} - \int \psi d \mu_{| \overline{V}} \int \varphi \frac{p_n^*( \lambda) }{\delta_l(p_n)} \right| =\left| \int \psi \Lambda_n(\varphi) d \mu_{| \overline{V}} - \int \psi d \mu_{| \overline{V}} \int \Lambda_n(\varphi) d \lambda \right| \\
&=\left| \int \psi \left(\Lambda_n(\varphi) - \langle \Lambda_n(\varphi) , \lambda \rangle \right) d \mu_{| \overline{V}}  \right| \leq   | \psi |_{L^s(\mu_{| \overline{V}})}  | \Lambda_n(\varphi) - \langle \Lambda_n(\varphi) , \lambda \rangle|_{L^r(\mu_{| \overline{V}})}\\
&\leq C_2  \left(    \frac{\delta_{l-1}(p_n)}{\delta_{l}(p_n)}  \right)^{\nu /2} | \psi |_{L^s(\mu_{| \overline{V}})} |\varphi |_{\mathcal{C}^{\nu}}.\\
\end{split}
\end{equation*}

This gives the second inequality in the Theorem. Notice that we can replace  $| \psi |_{L^s(\mu_{| \overline{V}})} $ with $| \psi |_{DSH}$ as in the previous remark.

\section{\bf Proof of Theorem \ref{Th2}}{\label{proof-Th2}}

\subsection{\bf Background in bifurcation theory:}{\label{background}}

In this paragraph, we follow the presentation of \cite{DGV} (see also \cite{DF} and \cite{Du}).

Let $\Lambda$ be a smooth complex quasi-projective variety and let $\widehat{f}: \Lambda \times \Pp^q \longrightarrow \Lambda \times \Pp^q$ be an algebraic family of endomorphisms of $\Pp^q$ of degree $d \geq 2$: $\widehat{f}$ is a morphism and $\widehat{f}(\lambda,z)=(\lambda, f_{\lambda}(z))$ where $f_{\lambda}$ is an endomorphism of $\Pp^q$ of algebraic degree $d$.

Assume that the family $\widehat{f}$ is endowed with $k$ marked points $a_1, \dots , a_k : \Lambda \longrightarrow \Pp^q$ where the $a_i$ are morphisms.

Let $\omega_{\Pp^q}$ be the Fubini-Study form on $\Pp^q$, $\pi_{\Lambda} : \Lambda \times \Pp^q \longrightarrow \Lambda$ and $\pi_{\Pp^q} : \Lambda \times \Pp^q \longrightarrow \Pp^q$ be the canonical projections.

If we denote $\widehat{\omega}:= (\pi_{\Pp^q})^* \omega_{\Pp^q}$, we have $\frac{\widehat{f}^* \widehat{\omega}}{d}= \widehat{\omega} + dd^c u$ with $u$ a smooth function (see \cite{DF} Proposition $3.1$). In the classical manner, the sequence 
$$\frac{(\widehat{f}^n)^* \widehat{\omega}}{d^n}= \widehat{\omega} + \sum_{i=0}^{n-1} dd^c \frac{u \circ \widehat{f}^i}{d^{i}}=\widehat{\omega} + dd^c u_n$$
(where $u_n= \sum_{i=0}^{n-1} \frac{u \circ \widehat{f}^i}{d^{i}}$) converges to a closed positive $(1,1)$-current $\widehat{T}= \widehat{\omega} + dd^c u_{\infty}$ on $\Lambda \times \Pp^q $ and this current has locally H\"older potential (see  \cite{DS08} Proposition $1.2.3$).

For $j=1, \cdots , k$, let $\Gamma_{a_j}$ be the graph of the marked point $a_j$ and we consider 

$$\mathfrak{a}=(a_1,\cdots,a_k): \Lambda \longrightarrow (\Pp^q)^k.$$

\begin{Def}

For $0 \leq j \leq k$, the bifurcation current $T_{a_j}$ of the point $a_j$ is the positive closed $(1,1)$-current on $\Lambda$ defined by

$$T_{a_j}=(\pi_{\Lambda})_* (\widehat{T} \wedge [ \Gamma_{a_j} ])$$

and we define the bifurcation current $T_{\mathfrak{a}}$ of the $k$-uple $\mathfrak{a}$ as

$$T_{\mathfrak{a}}=T_{a_1} + \cdot + T_{a_k}.$$

\end{Def}

For $n \in \Nn$, write 
$$\mathfrak{a}_n(\lambda)=(f_{\lambda}^n(a_1(\lambda)), \cdots , f_{\lambda}^n(a_k(\lambda)))=(a_{1,n}(\lambda) , \cdots ,  a_{k,n}(\lambda)) \mbox{    }   \mbox{   ,   } \mbox{    } \lambda \in \Lambda.$$

\begin{Lem} (See \cite{DF} Proposition 3.1 and \cite{DGV} Lemma 3.2).

For $1 \leq j \leq k$, the support of $T_{a_j}$ is the set of parameters $\lambda_0 \in \Lambda$ such that the sequence $(\lambda \longrightarrow f_{\lambda}^n(a_j(\lambda)))$ is not a normal family at  $\lambda_0$.

Moreover, there exists a locally uniformly bounded family of continuous functions $( u_{j,n})$ on $\Lambda$ such that

$$a_{j,n}^* (\omega_{\Pp^q})= d^n T_{a_j} + dd^c u_{j,n} \mbox{    }   \mbox{   on   } \mbox{    }  \Lambda.$$

\end{Lem}

To prove the last assertion, observe that 

$$a_{j,n}^* (\omega_{\Pp^q})={\pi_{\Lambda}}_* ((\widehat{f}^n)^* \widehat{\omega} \wedge \Gamma_{a_j})= {\pi_{\Lambda}}_* (d^n \widehat{T} \wedge  \Gamma_{a_j}) + {\pi_{\Lambda}}_*(dd^c (u_n - u_{\infty}) \wedge \Gamma_{a_j} d^n).$$

For $1 \leq j \leq k$ and $i \geq 1$, it follows that

$$a_{j,n}^* (\omega_{\Pp^q}^{i})= d^{ni} T_{a_j}^{i} + dd^c O(d^{(i-1)n})$$

on compact subset of $\Lambda$ and in particular $T_{a_j}^{q+1}=0$ on $\Lambda$ (see \cite{Ga} and \cite{Du} too).

Assume that $d_{\Lambda}:= \mbox{dim}(\Lambda)=qk$. Using the last property, the measure $T_{\mathfrak{a}}^{d_{\Lambda}}$ is equal to a constant multiplied by $T_{a_1}^{q} \wedge \cdots \wedge T_{a_k}^q$ and we define  

$$\mu_{\rm bif}:=T_{a_1}^{q} \wedge \cdots \wedge T_{a_k}^q.$$

This is the bifurcation measure of the $k$-uple $\mathfrak{a}=(a_1, \cdots, a_k)$. 

Let $\iota: \Lambda \hookrightarrow \Pp^m$ be an embedding of $\Lambda$ into a complex projective space. We identify $\Lambda$ with $\iota(\Lambda)$ and we denote by $\overline{\Lambda}$ the closure of $\Lambda$ in $\Pp^m$. Let $p_n: \overline{\Lambda} \longrightarrow (\Pp^q)^k $ be the meromorphic map obtained by taking the closure of the graph of $\mathfrak{a}_n$ in $\overline{\Lambda} \times (\Pp^q)^k$.

We can now proceed to the proof of Theorem \ref{Th2}, by using Theorem \ref{Th1} and pluripotential theory.

\subsection{\bf Proof of Theorem \ref{Th2}}

Pick $\psi \in DSH((\mathbb{P}^q)^k)$ and $\varphi \in C^{\nu}$ with compact support in $\overline{W}$.

Let $pr_j: (\mathbb{P}^q)^k \longrightarrow \Pp^q$ be the projection onto the $j$-th factor of the product $(\mathbb{P}^q)^k$ ($j=1, \cdots, k$) and consider $\Omega= \sum_{j=1}^{k} pr_j^* \omega_{\Pp^q}$. We will apply Theorem \ref{Th1} with $\lambda$ equal to $\Omega^{qk}$ normalized to be a probability, i.e., $\Omega^{qk}_{nor}= pr_1^* \omega_{\Pp^q}^q \wedge \cdots \wedge pr_k^* \omega_{\Pp^q}^q$ ($\omega_{\Pp^q}^j=0$ for $j>q$ and take the normalization $\int \omega_{\Pp^q}^q=1$).

In particular, we need estimates on $\delta_l(p_n)$ and $\delta_{l-1}(p_n)$ with $l=kq$, which are stated in the following Lemma:

\begin{Lem}{\label{Lem2}}
Suppose $\mu_{\rm bif} \neq 0$. Then there exists a positive constant $\epsilon$ such that

$$ \epsilon d^{kqn} \leq \delta_{kq}(p_n) \leq \frac{1}{\epsilon }d^{kqn} \mbox{    }  \mbox{  and  }  \mbox{    }  \delta_{kq-1}(p_n) \leq \frac{1}{\epsilon }d^{kqn-n}$$

for every $n \in \Nn$. Here, $\delta_{kq}(p_n):=   \int  p_n^*( \Omega^{qk})$ and $\delta_{kq-1}(p_n) := \int  p_n^*( \Omega^{qk-1}) \wedge \omega$, where $\omega$ is the Fubini-Study form on $\Pp^m(\Cc)$.

\end{Lem}

The proof of this Lemma is given at the end of the paragraph.

By using it with Theorem \ref{Th1}, there exists a positive constant $C$ (independent of $n \in \Nn$, $\psi \in DSH((\mathbb{P}^q)^k)$ and $\varphi \in C^{\nu}$) such that

$$\left| \int \psi(p_n) \varphi \frac{p_n^*( \mu_{| \overline{V}})}{d^{kqn}} - \int \psi d \mu_{| \overline{V}}  \int \varphi \frac{p_n^*( \Omega_{nor}^{qk}) }{d^{kqn}}   \right| \leq C d^{-n \nu/2}   | \psi |_{L^s(\mu_{| \overline{V}})} | \varphi |_{C^{\nu}},$$

so it remains to prove 

$$\left|  \int \psi d \mu_{| \overline{V}}  \int \varphi \frac{p_n^*( \Omega_{nor}^{qk}) }{d^{kqn}}  -  \int \psi d \mu_{| \overline{V}} \int \varphi  d \mu_{\rm bif}  \right| \leq C' d^{-n \nu/2}   | \psi |_{L^s(\mu_{| \overline{V}})} | \varphi |_{C^{\nu}},$$

for a positive constant $C'$ (independent of $n \in \Nn$, $\psi \in DSH((\mathbb{P}^q)^k)$ and $\varphi \in C^{\nu}$ with compact support in $\overline{W}$). We will again use interpolation theory between $\mathcal{C}^0$ and $\mathcal{C}^2$.

So, fix $\varphi \in \mathcal{C}^2(\overline{\Lambda})$ with compact support in $\overline{W}$. We have

\begin{equation*}
\begin{split}
 \int \varphi \frac{p_n^*( \Omega_{nor}^{qk}) }{d^{kqn}}&=\int_{\overline{W}}  \varphi \frac{p_n^*( \Omega_{nor}^{qk}) }{d^{kqn}}= \int_{\overline{W}} \varphi \frac{\mathfrak{a}_n^*(pr_1^* \omega_{\Pp^q}^q \wedge \cdots \wedge pr_k^* \omega_{\Pp^q}^q) }{d^{kqn}} \\
 &=\int_{\overline{W}} \varphi \frac{a_{1,n}^*  \omega_{\Pp^q}^q \wedge \cdots \wedge  a_{k,n}^*  \omega_{\Pp^q}^q }{d^{kqn}} \\
 &= \int \varphi \frac{ (d^n T_{a_1} + dd^c u_{1,n})^q \wedge \cdots \wedge (d^n T_{a_k} + dd^c u_{k,n})^q }{d^{kqn}} \\
\end{split}
\end{equation*}

where the $u_{j,n}$ are uniformly bounded on $\overline{W}$ (split $\overline{W}$ independently of $\varphi$ and use cut-off functions if necessary).

By Stokes' formula, the last integral is equal to

\begin{equation*}
\begin{split}
 & \int \varphi T_{a_1} \wedge (T_{a_1} + dd^c \frac{u_{1,n}}{d^n})^{q-1} \wedge \cdots \wedge (T_{a_k} + dd^c  \frac{u_{k,n}}{d^n})^q + \\
 &\int \frac{u_{1,n}}{d^n} dd^c \varphi  \wedge (T_{a_1} + dd^c \frac{u_{1,n}}{d^n})^{q-1} \wedge \cdots \wedge (T_{a_k} + dd^c  \frac{u_{k,n}}{d^n})^q= A + B\\
\end{split}
\end{equation*}

with obvious notations.

Let $0 \leq \theta_1\leq \cdots \leq \theta_{qk} \leq 1$ be smooth functions with compact support, $\theta_1 \equiv 1$ on a neighborhood of $\overline{W}$ and $\theta_{i+1} \equiv 1$ on a neighborhood of $\mbox{support}(\theta_{i})$ for $i=1, \cdots , qk-1$.

There exists a positive constant $C_1$ such that $-C_1 \theta_1 |\varphi|_{\mathcal{C}^2} \omega \leq u_{1,n} dd^c \varphi \leq C_1 \theta_1  |\varphi|_{\mathcal{C}^2} \omega$ for every $n$, where $\omega$ is the Fubini-Study form of $\Pp^m$. Hence

$$|B| \leq \frac{|\varphi|_{\mathcal{C}^2}}{d^n}   \int C_1 \theta_1 \omega  \wedge (T_{a_1} + dd^c \frac{u_{1,n}}{d^n})^{q-1} \wedge \cdots \wedge (T_{a_k} + dd^c  \frac{u_{k,n}}{d^n})^q.$$

Now, we prove that the previous integral is bounded by a constant independent of $n$. Indeed, write it as

\begin{equation*}
\begin{split}
 &C_1 \int \theta_1 \omega \wedge T_{a_1}  \wedge (T_{a_1} + dd^c \frac{u_{1,n}}{d^n})^{q-2} \wedge \cdots \wedge (T_{a_k} + dd^c  \frac{u_{k,n}}{d^n})^q +  \\
 &C_1 \int \frac{u_{1,n}}{d^n}  dd^c\theta_1 \wedge \omega   \wedge (T_{a_1} + dd^c \frac{u_{1,n}}{d^n})^{q-2} \wedge \cdots \wedge (T_{a_k} + dd^c  \frac{u_{k,n}}{d^n})^q .\\
\end{split}
\end{equation*}

As above, there exists a positive constant $C_2$ such that $-C_2 \theta_2  \omega \leq u_{1,n} dd^c \theta_1 \leq C_2 \theta_2  \omega$ and we iterate this process for both integrals, by using $\theta_1 , \cdots , \theta_{qk-1}$ successively.

At the end the integral $\int\theta_1 \omega  \wedge (T_{a_1} + dd^c \frac{u_{1,n}}{d^n})^{q-1} \wedge \cdots \wedge (T_{a_k} + dd^c  \frac{u_{k,n}}{d^n})^q$ is bounded above by a sum of terms like

$$C_1 \cdots C_l  \int \theta_l \omega^l  \wedge T_{a_1} ^{\alpha_1} \wedge \cdots \wedge T_{a_k}^{\alpha_k} $$

with $\alpha_1 + \cdots + \alpha_k = qk-l$ and $l=1, \cdots, qk$. All these integrals are bounded by a constant independent on $n$ since the potentials of the $T_{a_j}$ are continuous.

Hence there exists a positive constant $D$ such that 

$$|B| \leq \frac{D}{d^n} |\varphi|_{\mathcal{C}^2}.$$

Now for $A$ we follow the same method and we have

$$A= \int \varphi T_{a_1}^{q} \wedge \cdots \wedge T_{a_k}^q + \epsilon_n$$

with $| \epsilon_n| \leq \frac{D'}{d^n} |\varphi|_{\mathcal{C}^2}$.

We obtain

\begin{equation}{\label{eq}}
\left| \int \varphi \frac{p_n^*( \Omega_{nor}^{qk}) }{d^{kqn}} -  \int \varphi T_{a_1}^{q} \wedge \cdots \wedge T_{a_k}^q \right|  \leq \frac{D''}{d^n} |\varphi|_{\mathcal{C}^2}
\end{equation}

for $\mathcal{C}^2(\overline{\Lambda})$ maps $\varphi$ with compact support in $\overline{W}$.

When $\varphi$ is $\mathcal{C}^0(\overline{\Lambda})$, applying Lemma \ref{Lem2}, we have

$$\left| \int \varphi \frac{p_n^*( \Omega_{nor}^{qk}) }{d^{kqn}} -  \int \varphi T_{a_1}^{q} \wedge \cdots \wedge T_{a_k}^q \right|  \leq \left( \frac{c}{\epsilon}+  \mu_{\rm bif}(\Lambda) \right) |\varphi|_{\mathcal{C}^0}$$

where the constant $c > 0$ is such that  $ \Omega_{nor}^{qk} = c  \Omega^{qk}$.

Using interpolation theory for the linear operator

$$L: \{ \varphi \in \mathcal{C}^0  \mbox{     } \mbox{  with compact support in   } \mbox{     }   \overline{W} \}   \longrightarrow L^r(\mu_{| \overline{V}})$$

defined by 

$$L(\varphi)= \int \varphi \frac{p_n^*( \Omega_{nor}^{qk}) }{d^{kqn}} -  \int \varphi T_{a_1}^{q} \wedge \cdots \wedge T_{a_k}^q,$$

we obtain that there exists a constant $A_{\nu}$ such that

$$ \left| \int \varphi \frac{p_n^*( \Omega_{nor}^{qk}) }{d^{kqn}} -  \int \varphi T_{a_1}^{q} \wedge \cdots \wedge T_{a_k}^q   \right|_{L^r(\mu_{| \overline{V}})}$$ 

is bounded above by 

$$A_{\nu} \left( \left( \frac{c}{\epsilon}+  \mu_{\rm bif}(\Lambda) \right) \mu(\overline{V})^{1/r} \right)^{1- \nu/2} \left( \frac{D'' \mu(\overline{V})^{1/r}}{d^n} \right)^{\nu/2}  | \varphi |_{C^{\nu}},$$

for $\varphi \in C^{\nu}$ with compact support in $\overline{W}$.

Finally,

\begin{equation*}
\begin{split}
&\left|  \int \psi d \mu_{| \overline{V}}  \int \varphi \frac{p_n^*( \Omega_{nor}^{qk}) }{d^{kqn}}  -  \int \psi d \mu_{| \overline{V}} \int \varphi  d \mu_{\rm bif}  \right| 
= \left| \int \psi \left( \int \varphi \frac{p_n^*( \Omega_{nor}^{qk}) }{d^{kqn}} -  \int \varphi  d \mu_{\rm bif} \right)  d \mu_{| \overline{V}} \right| \\
&\leq  | \psi |_{L^s(\mu_{| \overline{V}})} \left| \int \varphi \frac{p_n^*( \Omega_{nor}^{qk}) }{d^{kqn}} -  \int \varphi T_{a_1}^{q} \wedge \cdots \wedge T_{a_k}^q \right|_{L^r(\mu_{| \overline{V}})} \leq C' d^{-n \nu/2} | \psi |_{L^s(\mu_{| \overline{V}})} | \varphi |_{C^{\nu}} \\
\end{split}
\end{equation*}

for $\varphi \in C^{\nu}$ with compact support in $\overline{W}$, and Theorem \ref{Th2} follows.

\vspace{0.5cm}

It now remains to prove Lemma \ref{Lem2}.

\vspace{0.5cm}

{\bf{Proof of Lemma \ref{Lem2}:}}

Since $\mu_{\rm bif} \neq 0$ by assumption, there exists a smooth function $0 \leq \theta_0 \leq 1$ with compact support in $\Lambda$ and $\int \theta_0 d \mu_{\rm bif} >0$.

Following the same method as in the previous proof with $\varphi=\theta_0$ and $W= \Lambda$, we obtain, as in inequality (\ref{eq}),

$$\left| \int \theta_0 \frac{p_n^*( \Omega_{nor}^{qk}) }{d^{kqn}} -  \int \theta_0 T_{a_1}^{q} \wedge \cdots \wedge T_{a_k}^q \right|  \leq \frac{C}{d^n}.$$

Therefore, using $ \Omega_{nor}^{qk} = c  \Omega^{qk}$, we have

\begin{equation*}
\begin{split}
\delta_{kq}(p_n)&:=   \int  p_n^*( \Omega^{qk}) = \frac{1}{c} \int  p_n^*( \Omega_{nor}^{qk}) \geq  \frac{1}{c}   \int \theta_0  p_n^*( \Omega_{nor}^{qk}) \\
& \geq \frac{ 1}{c} \left( \int \theta_0 T_{a_1}^{q} \wedge \cdots \wedge T_{a_k}^q - \frac{C}{d^n} \right)  d^{kqn} \geq \epsilon d^{kqn}\
\end{split}
\end{equation*}

where $\epsilon=\frac{1}{2c} \int \theta_0 T_{a_1}^{q} \wedge \cdots \wedge T_{a_k}^q $ and $n$ large enough (so for every $n$, up to reducing $\epsilon$).

To find the upper bounds of $\delta_{kq}(p_n)$ and $\delta_{kq-1}(p_n)$, we use Bezout's theorem like in \cite{S}.

First, $\delta_{kq}(p_n):=   \int  p_n^*( \Omega^{qk}) = \frac{1}{c} \int  p_n^*( \Omega_{nor}^{qk})$, so it is equal to the cardinal of $p_n^{-1}(b)$ with $b$ generic in $(\mathbb{P}^q)^k$ multiplied by the constant $\frac{1}{c}$. Thus, to compute it, we need to find the number of solutions to the equation $p_n(\lambda)=b$, or equivalently (by genericity), to the equation

 $$\mathfrak{a}_n(\lambda)=(f_{\lambda}^n(a_1(\lambda)), \cdots , f_{\lambda}^n(a_k(\lambda)))=b.$$

Let $b=(b_1, \cdots, b_k) \in (\mathbb{P}^q)^k$ be a generic point. For $j=1, \cdots, k$, write 

$$b_j=[b_{j,0}: \cdots : b_{j,q}]$$

 and

$$f_{\lambda}^n(a_j(\lambda))=[F_{\lambda ,0}^n(a_j(\lambda)): \cdots :  F_{\lambda ,q}^n(a_j(\lambda))]$$

where $F_{\lambda , 0}^n , \cdots , F_{\lambda ,q}^n$ are homogeneous polynomials of degree $d^n$ in $z_0, \cdots , z_q$ which define $f_{\lambda}^n$  (here  $[z_0: \cdots : z_q]$ are the coordinates in $\mathbb{P}^q$).

We are reduced to $k$ systems of equations as

\begin{equation}{\label{eq2}}
\left\{
\begin{array}{lll}
F_{\lambda ,0}^n(a_j(\lambda))= \frac{b_{j,0}}{b_{j,q}}  F_{\lambda, q}^n(a_j(\lambda)) \\
\vdots \\
F_{\lambda ,q-1}^n(a_j(\lambda))= \frac{b_{j,q-1}}{b_{j,q}}  F_{\lambda ,q}^n(a_j(\lambda)) \\
\end{array}
\right.
\end{equation}

(there is always a $b_{j,i_j} \neq 0$ and we assumed here $i_j=q$ to simplify the exposition).

The above equations are of degree $d_{j,0}(\lambda) d^n, \cdots , d_{j,q-1}(\lambda) d^n$ in $\lambda_0, \cdots , \lambda_m$, where $[\lambda_0: \cdots : \lambda_m]$ are coordinates in $\Pp^m$. The number of solutions to $\mathfrak{a}_n(\lambda)=b$ is finite since $p_n$ is dominant ($\int_{\Lambda} \frac{p_n^*( \Omega_{nor}^{qk}) }{d^{kqn}} >0$), so Bezout's inequality ($\overline{\Lambda}$ can have complete intersection or not) in $\Pp^m$ implies 

$$\delta_{kq}(p_n) \leq d(\lambda)d^{kqn} \leq \frac{1}{\epsilon} d^{kqn}$$

for every $n$, up to reducing $\epsilon$ if necessary.    

Now, we have to bound $\delta_{kq-1}(p_n):= \int  p_n^*( \Omega^{qk-1}) \wedge \omega= \int \mathfrak{a}_n ^*( \Omega^{qk-1}) \wedge \omega$ where $\omega$ is the Fubini-Study form on $\Pp^m(\Cc)$ (let us recall that $\iota: \Lambda \hookrightarrow \Pp^m$ and that we identify $\Lambda$ with $\iota(\Lambda)$).

Since $\Omega^{qk-1}= (\sum_{j=1}^{k} pr_j^* \omega_{\Pp^q})^{qk-1}=c(k,q) \sum_{j=1}^{k}  pr_1^* \omega_{\Pp^q}^q \wedge \cdots \wedge pr_j^* \omega_{\Pp^q}^{q-1} \wedge \cdots  \wedge pr_k^* \omega_{\Pp^q}^q$, we obtain

$$\delta_{kq-1}(p_n) =c(k,q) \sum_{j=1}^{k}  \int  a_{1,n}^*  \omega_{\Pp^q}^q \wedge \cdots \wedge a_{j,n}^*  \omega_{\Pp^q}^{q-1} \wedge  \cdots \wedge  a_{k,n}^*  \omega_{\Pp^q}^q \wedge \omega.$$

For cohomological reasons, to compute the integral above, we are left to bound the number of points in

$$\mathcal{E}= \overline{\Lambda} \cap a_{1,n}^{-1}(b_1) \cap \cdots \cap a_{j-1,n}^{-1}(b_{j-1}) \cap a_{j,n}^{-1}(L) \cap a_{j+1,n}^{-1}(b_{j+1}) \cap \cdots \cap a_{k,n}^{-1}(b_{k}) \cap H$$

where $b_1, \cdots , b_{j-1},  b_{j+1},  \cdots , b_{k}$ are generic points in $\Pp^q$, $L$ is a generic line in $\Pp^q$ and $H$ is a generic hyperplan in $\Pp^m$. Notice that this set is finite, since this is the intersection of a curve in $\overline{\Lambda}$ (as the preimage of a generic line of $(\Pp^q)^k$ by $p_n$, which is a dominant map) with a generic hyperplan $H$.

First, the points in $\mathcal{E}$ satisfy $k-1$ systems of equations as (\ref{eq2}). Then, the line $L$ is given by the intersection of $q-1$ hyperplans, so by $q-1$ equations of the type $\alpha_0^{i} z_0 + \cdots + \alpha_q^{i} z_q=0$ in $\Pp^q$ (for $i=1, \cdots , q-1$). Thus, the algebraic subset $a_{j,n}^{-1}(L)$  is given by the $q-1$ equations of the type

$$\alpha_0^{i} F_{\lambda ,0}^n(a_j(\lambda)) + \cdots + \alpha_q^{i} F_{\lambda ,q}^n(a_j(\lambda)) =0$$

which have degree $d'_{j,i}(\lambda) d^n$ (for $i=1, \cdots , q-1$).

In conclusion, by Bezout's inequality, the number of points in $\mathcal{E}$ is bounded above by $d'(\lambda) d^{(k-1)qn +(q-1)n}=d'(\lambda) d^{kqn-n}$, and the Lemma follows.

\bigskip

\noindent Henry De Th\'elin, Universit\'e Paris 13, Sorbonne Paris Nord, LAGA, CNRS (UMR 7539), F-93430, Villetaneuse, France.  

\noindent Email: {\tt dethelin@math.univ-paris13.fr}


\begin{thebibliography}{00}

\bibitem{BB1} G. Bassanelli and F. Berteloot, \textit{Bifurcation currents in holomorphic dynamics on $\Pp^k$}, J. Reine Angew. Math., {\bf 608} (2007), 201-235.

\bibitem{BD1} J.-Y. Briend and J. Duval, \textit{Exposants de Liapounoff
  et distribution des points p\'eriodiques d'un endomorphisme de $\Cc
  \Pp^k$}, Acta Math., {\bf 182} (1999), 143-157. 

\bibitem{BD2} J.-Y. Briend and J. Duval, \textit{Deux caract\'erisations
  de la mesure d'\'equilibre d'un endomorphisme de $\Pp^k(\Cc)$}, IHES
  Publ. Math., {\bf 93 } (2001), 145-159.

\bibitem{DeM} L. DeMarco,  \textit{Dynamics of rational maps: a current on the bifurcation locus}, Math. Res. Lett., {\bf 8} (2001), 57-66.

\bibitem{DGV} H. De Th\'elin, T. Gauthier and G. Vigny, \textit{The bifurcation measure has maximal entropy}, Isr. J. Math., {\bf 235} (2020), 213-243.

\bibitem{DGV2} H. De Th\'elin, T. Gauthier and G. Vigny, \textit{Parametric Lyapunov exponents}, Bull. Lond. Math. Soc., {\bf 53} (2021), 660-672.

\bibitem{DS06} T.-C. Dinh and N. Sibony, \textit{Distribution des valeurs de transformations m\'eromorphes et applications}, Comment. Math. Helv., {\bf 81} (2006), 221-258.

\bibitem{DS08} T.-C. Dinh and N. Sibony, \textit{Dynamics in several complex variables: endomorphisms of projective spaces and polynomial-like mappings}, Lecture Notes in Mathematics, {\bf 1998} (2010), 165-294.

\bibitem{DSN} T.-C. Dinh, V.-A. Nguyen and N. Sibony, \textit{Exponential estimates for plurisubharmonic functions}, J. Diff. Geom., {\bf 84} (2010), 465-488.

\bibitem{Du} R. Dujardin, \textit{The supports of higher bifurcation currents}, Ann. Fac. Sci. Toulouse Math., {\bf 22} (2013), 445-464.

\bibitem{DF} R. Dujardin and C. Favre, \textit{Distribution of rational maps with a preperiodic critical point}, Am. J. Math., {\bf 130} (2008), 979-1032.

\bibitem{FS} J.E. Forn{\ae}ss and N. Sibony, \textit{Complex dynamics in higher dimensions}, Complex Potential Theory (Montreal, PQ, 1993), NATO Adv. Sci. Inst. Ser. C Math. Phys. Sci., {\bf 439}, Kluwer, Dordrecht (1994), 131-186.

\bibitem{FS1} J.E. Forn{\ae}ss and N. Sibony, \textit{Complex dynamics in higher dimension I}, Ast\'erisque {\bf 222} (1994), 201-231.

\bibitem{Ga} T. Gauthier, \textit{Strong bifurcation loci of full Hausdorff dimension}, Ann. Sci. \'Ec. Norm. Sup\'er., {\bf 45} (2012), 947-984.

\bibitem{GKN} D. Ghioca, H. Krieger and K. Nguyen, \textit{A case of the dynamical Andr\'e-Oort conjecture}, Int. Math. Res. Not., {\bf 3} (2016), 738-758.

\bibitem{GS} J. Graczyk and G. \'Swiatek, \textit{Lyapunov exponent and harmonic measure on the boundary of the connectedness locus}, Int. Math. Res. Not., {\bf 16} (2015), 7357-7364.

\bibitem{L} Y. Luo, \textit{On the inhomogeneity of the Mandelbrot set}, Int. Math. Res. Not., {\bf 8} (2021), 6051-6076.

\bibitem{M} M. M\'eo, \textit{Inverse image of a closed positive current by a surjective analytic map}, C. R. Acad. Sci. Paris, {\bf 12} (1996), 1141-1144.

\bibitem{S} N. Sibony, \textit{Dynamique des applications rationnelles de $\Pp^k$}, Panor. Synth\`{e}ses, {\bf 8} (1999), 97-185.




\end{thebibliography}
\end{document}